\documentclass[number,citesort,MSNbibl,dvips]{arxbj}
\usepackage{graphicx}

\aid{0}
\volume{18}
\issue{2}
\pubyear{2012}
\firstpage{455}
\lastpage{475}
\doi{10.3150/10-BEJ343}

\makeatletter

\def\bsuffix #1{#1}

\newcommand{\eqref}[1]{(\ref{#1})}

\newtheorem{lemma}{Lemma}[section]
\newtheorem{theorem}[lemma]{Theorem}
\newtheorem{cor}[lemma]{Corollary}
\newtheorem{prop}[lemma]{Proposition}
\newremark{exam}[lemma]{Example}

\newcommand{\R}{\mathbb{R}}
\newcommand{\N}{\mathbb{N}}
\newcommand{\norm}[1]{\Vert #1\Vert}
\newcommand{\abs}[1]{\vert #1\vert}
\newcommand{\set}[1]{\{#1\}}

\newcommand{\bfzero}{\mathbf{0}}
\newcommand{\bfone}{\mathbf{1}}
\newcommand{\bfa}{\mathbf{a}}
\newcommand{\bfb}{\mathbf{b}}
\newcommand{\bfc}{\mathbf{c}}
\newcommand{\bfd}{\mathbf{d}}
\newcommand{\bfQ}{\mathbf{Q}}
\newcommand{\bfS}{\mathbf{S}}
\newcommand{\bfT}{\mathbf{T}}
\newcommand{\bfu}{\mathbf{u}}
\newcommand{\bfV}{\mathbf{V}}
\newcommand{\bfX}{\mathbf{X}}
\newcommand{\bfY}{\mathbf{Y}}
\newcommand{\bfZ}{\mathbf{Z}}
\newcommand{\bfSigma}{\bolds{\Sigma}}

\makeatother

\begin{document}
\begin{frontmatter}

\title{A multivariate piecing-together approach with an application to operational loss data}
\runtitle{Multivariate piecing-together}

\begin{aug}
\author[1]{\fnms{Stefan} \snm{Aulbach}\thanksref{1,e1}\ead[label=e1,mark]{stefan.aulbach@uni-wuerzburg.de}},
\author[2]{\fnms{Verena} \snm{Bayer}\thanksref{2}\ead[label=e2]{verena.bayer@uni-wuerzburg.de}}
\and
\author[1]{\fnms{Michael} \snm{Falk}\corref{}\thanksref{1,e3}\ead[label=e3,mark]{michael.falk@uni-wuerzburg.de}}
\runauthor{S. Aulbach, V. Bayer and M. Falk}
\address[1]{Institute of Mathematics, Emil-Fischer-Str. 30, D-97074 W\"{u}rzburg, Germany.\\
\printead{e1,e3}}
\address[2]{Chair of Econometrics, Sanderring 2, D-97070 W\"{u}rzburg, Germany.\\
\printead{e2}}
\end{aug}

\received{\smonth{9} \syear{2009}}
\revised{\smonth{8} \syear{2010}}

%
\begin{abstract}
The univariate piecing-together approach (PT) fits a univariate
generalized Pareto distribution (GPD) to the upper tail of a given
distribution function in a continuous manner. We propose a~multivariate
extension. First it is shown that an arbitrary copula is in the domain
of attraction of a multivariate extreme value distribution if and only
if its upper tail can be approximated by the upper tail of a~multivariate
GPD with uniform margins.

The multivariate PT then consists of two
steps: The upper tail of a given copula $C$ is cut off and substituted
by a multivariate GPD copula in a continuous manner.
The result is again a~copula.
The other step consists of the transformation of each margin of this
new copula by a~given univariate distribution function.

This provides, altogether, a multivariate distribution function
with prescribed margins whose copula coincides in its central part
with $C$ and in its upper tail with a GPD copula.

When applied to data, this approach also enables the
evaluation of a wide range of rational scenarios for the
upper tail of the underlying distribution function in
the multivariate case. We apply this approach to operational loss data
in order to
evaluate the range of operational risk.
\end{abstract}

%

\end{frontmatter}

\section{Introduction}

The peaks over threshold approach (POT) shows that the upper tail of a
univariate distribution function $F$ can reasonably be approximated
only by that of a generalized Pareto distribution (GPD). This result
goes back to Balkema and de Haan \cite{balkemahaan74} and Pickands~\cite{pickands75}.
A univariate GPD $W$ is derived from an extreme value distribution (EVD)
$G$ by the equality
\[
W(x)=1+\log(G(x)),\qquad 1/\mathrm{e}\le G(x),
\]
where, with a shape parameter $\alpha>0$, the family of standardized
EVD is given by
%
\begin{eqnarray}\label{eqnunivariateEVD}
G_{1,\alpha}(x)&=&\exp(-x^{-\alpha}),\qquad  x>0,\nonumber\\
G_{2,\alpha}(x)&=&\exp(-(-x)^{\alpha}),\qquad  x\le0,\\
G_3(x)&=&\exp(-\mathrm{e}^{-x}),\qquad  x\in\R,\nonumber
\end{eqnarray}
being the Fr\'{e}chet, (reverse) Weibull and Gumbel case of an EVD.

The family of univariate standardized GPD is, consequently, given by
\begin{eqnarray*}
W_{1,\alpha}(x)&=&1-x^{-\alpha},\qquad  x\ge1,\\
W_{2,\alpha}(x)&=&1-(-x)^{\alpha},\qquad -1\le x\le0,\\
W_3(x)&=&1-\exp(-x),\qquad  x\ge0,
\end{eqnarray*}
being the Pareto, beta and exponential case of a GPD.

If $X$ is a univariate random variable with distribution function $F$,
then the distribution function $F^{[x_0]}$ of $X$, conditional on the
event $X>x_0$, is given by
\begin{eqnarray*}
F^{[x_0]}(x)&=&P(X\le x\mid X> x_0)\\
&=&\frac{F(x)-F(x_0)}{1-F(x_0)},\qquad  x\ge x_0,
\end{eqnarray*}
where we require $F(x_0)< 1$. The POT approach shows that $F^{[x_0]}$
can reasonably be approximated only by a GPD with appropriate shape,
location and scale parameter~$W_{\gamma,\mu,\sigma}$. Note that
\begin{eqnarray*}
F(x)&=&\bigl(1-F(x_0)\bigr)F^{[x_0]}(x)+F(x_0)\\
&\approx&\bigl(1-F(x_0)\bigr)W_{\gamma,\mu,\sigma}(x)+F(x_0),\qquad  x\ge x_0.
\end{eqnarray*}
The piecing together approach (PT) now consists in replacing the
distribution function~$F$ by
%
\begin{equation}\label{eqnunivariatept}
F_{x_0}^*(x)=
\cases{ F(x),&\quad $x< x_0$,\cr
\bigl(1-F(x_0)\bigr)W_{\gamma,\mu,\sigma}(x)+F(x_0),&\quad $x\ge x_0$,
}
\end{equation}
where the shape, location and scale parameters $\gamma$, $\mu$,
$\sigma$ of the GPD are typically estimated from given data. This
modification aims at a more precise investigation of the upper end of
the data.

Replacing $F$ in \eqref{eqnunivariatept} by the empirical distribution
function $\hat F_n$ of $n$ independent copies of $X$ offers in
particular a semi-parametric approach to the estimation of high
quantiles $F^{-1}(q)=\inf\set{t\in\R\dvt F(t)\ge q}$ outside the range of
given data; see, for example, Section~2.3 of Reiss and Thomas \cite{reth}.

In this paper we propose an extension of the PT in
\eqref{eqnunivariatept} to higher dimensions. When applied to data,
this approach also enables the evaluation of a wide range of rational
scenarios for the upper tail of the underlying distribution function in
the multivariate case. This will be exemplified in Section
\ref{secoperationalrisk} for operational loss data, where we simulate
different scenarios for risk parameters such as the value at risk or
the expected shortfall. In Section \ref{secmultivariategpd} we provide
the basic mathematics for our PT approach. We will show that an
arbitrary copula can reasonably be approximated in its upper tail only
by a GPD with uniform margins.

The multivariate PT approach, which will be established in Section
\ref{secmultivariatepiecingtogether}, now consists of two steps:
\begin{longlist}
\item The upper tail of a given $m$-dimensional copula $C$ is cut off
and substituted by the upper tail of multivariate GPD copula in a
continuous manner such that the result is again a copula.
\item The
other step consists of the transformation of each margin of this new
copula by a given univariate distribution function $F_i^*$, $1\le i\le
m$.
\end{longlist}
This provides, altogether, a multivariate distribution function with
prescribed margins~$F_i^*$, whose copula coincides in its central part
with $C$ and in its upper tail with a GPD copula. In Section~\ref{secsimulation} we will simulate the effects that the combination
of univariate and multivariate PT has on quantile functions and mean
excess functions or, in terms of risk analysis, on value at risk and
expected shortfall. It turns out that in our specific model, which will
be specified in Section~\ref{secoperationalrisk}, the application of
the multivariate PT approach leads to a rising expected shortfall while
the value at risk keeps (up to a~level of 99.9\%) nearly unchanged.

Instead of fitting a GPD to the upper tail of a distribution,
estimation of rare events in the multivariate case can also be based on
the fact that the exponent measure pertaining to a multivariate GPD is
homogeneous; see de Haan and Sinha \cite{hasi99} and de Haan and Ronde
\cite{haro98}
for details.

For recent accounts of basic and advanced topics of extreme value
theory and statistics, see the monographs by Reiss and Thomas \cite{reth},
de Haan and Ferreira~\cite{hafe06} and
Resnick~\cite{resnick06}.\vspace*{-1pt}

\section{Multivariate GPD}\vspace*{-1pt}\label{secmultivariategpd}

In this section we provide the mathematics underlying our PT approach,
which will be established in Section
\ref{secmultivariatepiecingtogether}.

Let $F$ be an arbitrary $m$-dimensional distribution function that is
in the domain of attraction of an $m$-dimensional EVD $G$; that is,
there exist norming constants $\bfa_n>\bfzero$, $\bfb_n\in\R^m$ such
that
%
\begin{equation}\label{eqnmultivariatedomainofattraction}
F^n(\bfa_n\mathbf{x}+\bfb_n)\mathop{\to}\limits_{n\to\infty}G(\mathbf{x}),\qquad  \mathbf{x}\in\R^m,
\end{equation}
where all operations on vectors are meant componentwise. The
distribution function $G$ is max stable; that is, there exist norming
constants $\bfc_n>\bfzero$, $\bfd_n\in\R^m$ with
\[
G^n(\bfc_n\mathbf{x}+\bfd_n)=G(\mathbf{x}),\qquad \mathbf{x}\in\R^m.
\]
The one-dimensional margins $G_i$ of $G$ are up to scale and location
parameters univariate EVD in \eqref{eqnunivariateEVD}.\vadjust{\goodbreak}

It is well known that \eqref{eqnmultivariatedomainofattraction} is
equivalent with convergence of the univariate margins \textit{together}
with convergence of the \textit{copulas}
%
\begin{equation}\label{eqnconvergenceofcopulas}
\lim_{n\to\infty} C_F^n(\bfu^{1/n})= C_G(\bfu) =G(G_1^{-1}(u_1), \dots,
G_m^{-1}(u_m)),\qquad \bfu\in(0,1)^m,\vspace*{-1pt}
\end{equation}
(Deheuvels \cite{deheuvels78,deheuvels84}, Galambos \cite{ga87}). For a recent account on
copulas, see Nelsen \cite{nelsen06}. Some elementary computations as in Falk
\cite{falk08}, Section 6, or de Haan and Ronde \cite{haro98}, Section 4.2, entail that
convergence \eqref{eqnconvergenceofcopulas} is equivalent with
%
\begin{equation}\label{eqnconvergencetostabletaildepfunc}
\lim_{t\downarrow0} \frac1t\bigl(1-C_F(\mathbf{1}+t\mathbf{x})\bigr)=
l_G(\mathbf{x}):=-\log(C_G(\exp(\mathbf{x}))),\qquad \mathbf{x}\le\bfzero,\vspace*{-1pt}
\end{equation}
where $l_G$ is known as the \textit{stable tail dependence function}
introduced by Huang \cite{huang92}. For a~detailed discussion of the stable
tail dependence function, see Beirlant \textit{et al.} \cite{begosete04}. The stable tail
dependence function is homogeneous $tl_G(\mathbf{x})=l_G(t\mathbf{x})$, $t\ge0$,
and, thus, \eqref{eqnconvergencetostabletaildepfunc}
becomes
\[
\frac{1-C_F(\mathbf{1}+t\mathbf{x})-l_G(t\mathbf{x})}t\to_{t\downarrow0}0.\vspace*{-1pt}
\]

Observe that $l_G(\mathbf{x})=1-H(\mathbf{x})$, $\mathbf{x}\le\bfzero$, where
$H$ is a
multivariate GP \textit{function} with uniform margins $H_i(x)=1+x$,
$x\le0$, $i\le m$; that is,
\[
H(\mathbf{x})=1+\log(\widetilde G(\mathbf{x})),\qquad \mathbf{x}\le\bfzero,\vspace*{-1pt}
\]
and $\widetilde G$ is a multivariate EVD with negative exponential
margins $\widetilde G_i(x)=\exp(x)$, $x\le0$, $i\le m$.

We call in general an $m$-dimensional \textit{distribution function}
$W$ a multivariate GPD if its upper tail coincides with a GP function;
that is, there exist a multivariate EVD $G$ and a vector
$\mathbf{x}_0\in\R^m$ with $G(\mathbf{x}_0)<1$ such that
%
\begin{equation}\label{eqnmultivariateGPD}
W(\mathbf{x})=1+\log(G(\mathbf{x})),\qquad \mathbf{x}\ge\mathbf{x}_0.\vspace*{-1pt}
\end{equation}

Note that $H(\mathbf{x})=1+\log(G(\mathbf{x}))$, $G(\mathbf{x})\ge1/\mathrm{e}$, does
not define a
distribution function unless $m\in\set{1,2}$; see Michel \cite{michel08}, Theorem
6. We, therefore, call $H$ a GP function. It is, actually, a
quasi-copula (Alsina \textit{et al.} \cite{alnesch93}, Genest \textit{et al.} \cite{gemorose99} and
Section 5.1 in Falk \textit{et al.} \cite{fahure04}). Lemma 5.1.5 in Falk \textit{et al.}
\cite{fahure04} implies, on the other hand, that for any GP function there
exists a distribution function $W$ satisfying
\eqref{eqnmultivariateGPD}.

The preceding considerations together with elementary computations
entail now the following characterization of domains of attraction in
terms of a GPD. By $\norm{\cdot}$ we denote an arbitrary norm on
$\R^m$.\vspace*{-2pt}

\begin{theorem}\label{theocharacterizationofdomain}
An arbitrary distribution function $F$ is in the domain of attraction
of a~multivariate EVD $G$ if and only if this is true for the
univariate margins and if there exists a GPD $W$ with ultimately
uniform margins $W_i(x)=1+x$, $x_0\le x\le0$, $i\le m$, such that
\[
C_F(\mathbf{y})=W(\mathbf{y}-\mathbf{1})+\mathrm{o}(\norm{\mathbf{y}-\mathbf{1}})\vspace*{-1pt}
\]
uniformly for $\mathbf{y}\in[0,1]^m$.\vspace*{-2pt}\vadjust{\goodbreak}
\end{theorem}

We have the following equivalences for an arbitrary copula $C$ to lie
in the domain of attraction of an EVD.

\begin{cor}\label{coroequivalencesofdomain}
$C$ is in the domain of attraction of an EVD $G$
\begin{itemize}[$\iff$]
\item[$\iff$] There exists a GPD $W$ with ultimately uniform margins
such that
\[
C(\mathbf{y})=W(\mathbf{y}-\mathbf{1})+\mathrm{o}( \norm{\mathbf{y}-\mathbf{1}}),
\]
uniformly for $\mathbf{y}\in[0,1]^m$. In this case $W(\mathbf{x})=1+\log
(G(\mathbf{x}))$,
$\mathbf{x}_0\le\mathbf{x}\le\bfzero\in\R^m$.
\item[$\iff$] There exists a norm
$\norm{\cdot}_D$ on $\R^m$ such that
\[
C(\mathbf{y})=1-\norm{\mathbf{y}-\mathbf{1}}_D+ \mathrm{o}(\norm{\mathbf{y}-\mathbf{1}}_D),
\]
uniformly for $\mathbf{y}\in[0,1]^m$. In this case
$G(\mathbf{x})=\exp(-\norm{\mathbf{x}}_D)$, $\mathbf{x}\le\bfzero$.
\end{itemize}
\end{cor}

Recall that all norms on $\R^m$ are equivalent and, thus,
$\mathrm{o}(\norm{\mathbf{y}-\mathbf{1}}_D)$ in the second equivalence above can be
substituted by $\mathrm{o}(\norm{\mathbf{y}-\mathbf{1}})$ with an arbitrary norm
$\norm\cdot$ on $\R^m$.

The preceding results show that the upper tail of the copula $C_F$ of a
distribution function $F$ can reasonably be approximated only by that
of a GPD $W$ with ultimately uniform margins. To the best of our
knowledge, this provides new insight into the significance of
multivariate GPD. But it is in accordance with Rootz\'{e}n and Tajvidi
\cite{rootzen05}, who showed that, in the multivariate case, modelling
exceedances of a random variable over a high threshold can rationally
be done only by a multivariate GPD.

\begin{pf*}{Proof of Corollary \ref{coroequivalencesofdomain}}
It is well known that a GPD $W$ with ultimately uniform margins can be
written as
\[
W(\mathbf{x})=1-\norm{\mathbf{x}}_D,\qquad \mathbf{x}_0\le\mathbf{x}\le\bfzero,
\]
where $\norm{\cdot}_D$ is a norm on $\R^m$ with particular properties,
called a $D$-norm; see Section~4.4 in Falk \textit{et al.} \cite{fahure04}. In
particular, $G(\mathbf{x})= \exp(-\norm{\mathbf{x}}_D)$, $\mathbf{x}\le
\bfzero$, defines
an EVD on $\R^m$. If $C(\mathbf{y})=W(\mathbf{y}-\mathbf{1})+\mathrm{o}( \norm{\mathbf
{y}-\mathbf{1}})$,
$\mathbf{y}\in[0,1]^m$, for some norm $\norm{\cdot}$ on $\R^m$, then
\begin{eqnarray*}
C^n\biggl(\mathbf{1}+\frac yn\biggr)&=&\biggl(1-\frac1n \norm{\mathbf{y}}_D+\mathrm{o}\biggl(\frac1n \norm
{\mathbf{y}}\biggr)\biggr)^n\\
&\mathop{\to}\limits_{n\to\infty}& \exp(-\norm{\mathbf{y}}_D)=G(\mathbf{y}),\qquad \mathbf {y}\le\bfzero.
\end{eqnarray*}
Together with Theorem \ref{theocharacterizationofdomain} this implies
Corollary \ref{coroequivalencesofdomain}.
\end{pf*}

In the final equivalence of Corollary \ref{coroequivalencesofdomain},
the norm can obviously be computed as
\[
\norm{\mathbf{x}}_D=\lim_{t\downarrow0} \frac{1-
C(\mathbf{1}+t\mathbf{x})}t=l(\mathbf{x}),\qquad \mathbf{x}\le\bfzero;
\]
that is, it is the stable tail dependence function. It turns out that
any stable tail dependence function is actually a norm. This explains
why it is a convex function and homogeneous of order one.\vadjust{\goodbreak}

\begin{exam}
Take an arbitrary Archimedean copula
\[
C_\varphi(\bfu)=\varphi^{-1}\bigl(\varphi(u_1)+
\cdots+\varphi(u_m)\bigr),\vspace*{-2pt}
\]
where the generator $\varphi\dvtx(0,\infty)\to[0,\infty)$ is a continuous
function that is strictly decreasing on $(0,1]$, $\varphi(1)=0$,
$\lim_{x\downarrow0}\varphi(x)=\infty$ and
$\varphi^{-1}(t)=\inf\set{x>0\dvt\varphi(x)\le t}$, $t\ge0$.

Note that $C_\varphi$ is not automatically a copula for each function
$\varphi\dvt(0,\infty)\to[0,\infty)$ as above. While, in the bivariate
case $m=2$, convexity of $\varphi^{-1}$ is a necessary and sufficient
condition, this is no longer true in higher dimension $m\ge3$.
Instead, $C_\varphi$ is, for general dimension $m\ge2,$ a copula if
and only if $\varphi^{-1}$ is differentiable up to order $m-2$, the
derivatives satisfy $(-1)^k (\varphi^{-1})^{(k)}(x)\ge0$,
$k=0,\dots,m-2$, $x\in(0,\infty)$ and further if $(-1)^{m-2}
(\varphi^{-1})^{(m-2)}$ is non-increasing and convex in $(0,\infty)$;
see McNeil and Ne\v{s}lehov\'{a}~\cite{mcneilneslehova09}, Theorem 2.2.

If $\varphi$ is differentiable from the left in $x=1$ with left
derivative $\varphi'(1-)\not=0$, then
\[
\lim_{t\downarrow0}\frac{1-C_\varphi(\mathbf{1}+t\mathbf{x})}t=\sum _{i\le
m}\abs{x_i}=\norm{\mathbf{x}}_1,\qquad \mathbf{x}\le\bfzero;\vspace*{-2pt}
\]
that is, each Archimedean copula with a generator $\varphi$ as above is
in the domain of attraction of the EVD $G(\mathbf{x})=\exp(-\norm{\mathbf{x}}_1)$,
$\mathbf{x}\le\bfzero$, with independent margins. The margins of~$C_\varphi$
are, therefore, tail independent; that is, the tail dependence
parameters vanish:
\[
\chi(i,j):=\lim_{x\uparrow1}P(U_i>x\mid U_j>x)=0,\qquad 1\le i\not =j\le m,\vspace*{-2pt}
\]
where the random vector $(U_1,\dots,U_m)$ follows the distribution
function $C_\varphi$. For a discussion of the tail dependence
parameter and further literature, see Section 6.1 in Falk \textit{et al.}
\cite{fahure04}.\vspace*{-3pt}
\end{exam}

The preceding considerations concern, for example, the Clayton and the
Frank copula, which have generators $\varphi_C(t) =\vartheta^{-1}
(t^{-\vartheta}-1)$ and $\varphi_F(t)=-\log ((\exp(-\vartheta
t)-1)/\allowbreak
(\exp(-\vartheta)-1) )$, $\vartheta>0$, but not the Gumbel copula with
parameter $\lambda>1$, which has generator $\varphi
_G(t)=-(\log(t))^\lambda$, $\lambda\ge1$, $0<t\le1$.

Any multivariate EVD $G$ has univariate EVD margins and any
multivariate GPD $W$ has univariate GPD margins in its upper tail. We
can transform an arbitrary multivariate EVD to an EVD with negative
exponential margins by just transforming the margins. Equally, we can
transform an arbitrary $W$ to a GPD with uniform margins by just
transforming the margins. This transformation can also be done
backwards; see Section~5.6 in Falk \textit{et al.} \cite{fahure04}. We will, therefore,
consider in what follows multivariate GPD derived from an EVD $G$ with
negative exponential margins. For a recent account on multivariate GPD,
see Michel~\cite{michel06}.

From the de Haan--Resnick--Pickands representation of a multivariate
EVD, it is well known that a function $G$ on $(-\infty,0]^m$ is the
distribution function of an EVD with negative standard exponential
margins $G_i(x)=\exp(x)$, $x\le0$, $i\le m$, if and only if it can be
represented as
\[\label{eqndehaanresnickpickandsrepresentation}
G(\mathbf{x})=\exp\biggl(\int_{S_m}\min_{i\le m}(x_it_i)\mu(\mathrm{d}\mathbf{t})\biggr),\qquad \mathbf{x}\le
\bfzero,\vspace*{-2pt}\vadjust{\goodbreak}
\]
where $\mu$ is a finite measure on $S_m:=\set{\mathbf{t}\ge\bfzero
\dvt\sum_{i\le
m}t_i=1}$, called an \textit{angular measure}, with the characteristic
property $\int_{S_m}t_i\mu(\mathrm{d}\mathbf{t})=1$, $i\le m$; see Section 4.2 in Falk
\textit{et al.} \cite{fahure04}. Note that this integrability condition on $\mu$
implies that $\mu(S_m)=\int_{S_m}1\,\mathrm{d}\mu=\int_{S_m}\sum_{i\le
m}t_i\mu(\mathrm{d}\mathbf{t})=\sum_{i\le
m}\int_{S_m}t_i\mu(\mathrm{d}\mathbf{t})=m$.\vspace*{2pt}

As a consequence we obtain that a multivariate GPD $W$ with standard
uniform margins $1-W_i(x)=x$, $i\le m$,
in a left neighborhood of $\bfzero\in\R^m$
can be represented as
\begin{eqnarray}\label{eqndehaanresnickpickandsrepresentationofgpd}
W(\mathbf{x})&=&1+\biggl(\sum_{j\le m}x_j\biggr)\int_{S_m}\max_{i\le m}(\tilde
x_it_i)\mu(\mathrm{d}\mathbf{t})\nonumber\\ [-9pt]\\ [-9pt]
&=:&1+ \biggl(\sum_{j\le m}x_j\biggr)D(\tilde x_1,\dots, \tilde x_{m-1})\nonumber\vspace*{-2pt}
\end{eqnarray}
for $\mathbf{x}_0\,{\le}\,\mathbf{x}\,{\le}\,\bfzero$, where $\mu$ is as above,
$\tilde
x_i\,{=}\,x_i/\sum_{j\le m}x_j$ and $D\dvt\{\bfu\,{\in}\,[0,1]^{m-1}\dvt\sum_{j\le
m-1}u_j\,{\le}\allowbreak1\}\to[1/m,1]$ is a \textit{Pickands dependence function}
(Section 4.3 in Falk \textit{et al.} \cite{fahure04}).

The following result characterizes a GPD with uniform margins in terms
of random variables. It provides an easy way to generate a multivariate
GPD, thus extending the bivariate approach proposed by
Buishand \textit{et al.} \cite{buhazh06} to an arbitrary dimension. Recall that an arbitrary
multivariate GPD can be obtained from a GPD with ultimately uniform
margins by just transforming the margins. For a recent account on
simulation techniques of multivariate GPD, see Michel
\cite{michel07}.\vspace*{-3pt}

\begin{prop}\label{propgenerationofGPD}
\begin{longlist}
\item Let $W$ be a multivariate GPD with standard uniform margins in a
left neighborhood of $\bfzero\in\R^m$. Then there is a random vector
$\bfZ=(Z_1,\dots,Z_m)$ with $Z_i\in[0,m]$ and $E(Z_i)=1$, $i\le m$, and
a vector $(-1/m,\dots,-1/m)\le\mathbf{x}_0<\bfzero$ such that
\[
W(\mathbf{x})=P\biggl(-U\biggl(\frac1{Z_1},\dots,\frac1{Z_m}\biggr)\le\mathbf{x}\biggr),\qquad
\mathbf{x}_0\le\mathbf{x}\le\bfzero,\vspace*{-2pt}
\]
where the random variable $U$ is uniformly distributed on $(0,1)$ and
independent of $\bfZ$.
\item\label{itemdefinitionofaGPDbyarv}The
random
vector $-U(1/Z_1,\dots,1/Z_m)$ follows a GPD with
standard uniform margins in a left neighborhood of
$\bfzero\in\R^m$ if $U$ is independent of
$\bfZ=(Z_1,\dots,Z_m)$ and $0\le Z_i\le c_i$ a.s. with
$E(Z_i)=1$, $i\le m$, for some $c_1,\dots,c_m\ge1$.\vspace*{-3pt}
\end{longlist}
\end{prop}

Note that the case of a GPD $W$ with arbitrary uniform margins
$W_i(x)=1-a_ix$ in a~left neighborhood of $\bfzero$ with arbitrary
scaling factors $a_i>0$, $i\le m$, immediately follows from the
preceding result by substituting $Z_i$ by $a_iZ_i$.\vspace*{-3pt}

\begin{pf*}{Proof of Proposition \ref{propgenerationofGPD}}
First we establish part (i). From representation
\eqref{eqndehaanresnickpickandsrepresentationofgpd} we obtain that for
$\mathbf{x}$ in a left neighborhood of $\bfzero\in\R^m$
\[
W(\mathbf{x})=1+\biggl(\sum_{j\le m}x_j\biggr)\int_{S_m}\max_{i\le m}(\tilde
x_it_i)\mu(\mathrm{d}\mathbf{t})\vspace*{-2pt}
\]
with some measure $\mu$ on $S_m$ such that $\mu(S_m)=m$ and
$\int_{S_m}t_i\mu(\mathrm{d}\mathbf{t})=1$, $i\le m$.\vadjust{\goodbreak}

Now $\tilde\mu(\cdot)=\mu(\cdot)/m$ defines a probability measure on
$S_m$. Let $\bfT=(T_1,\dots,T_m)$ be a random vector with values in
$S_m$ that has distribution $\tilde\mu$ and put $\bfZ:=m\bfT$. Then
$\bfZ\in[0,m]^m$ and $E(Z_i)=\int_{S_m}t_i\mu(\mathrm{d}\mathbf{t})=1$, $i\le
m$. We
have, further, for $\mathbf{x}\le\bfzero\in\R^m$ with $x_j\ge-1/m$,
$j\le
m$,
\begin{eqnarray*}
&&P\biggl(-U\biggl(\frac1{Z_1},\dots,\frac1{Z_m}\biggr)\le\mathbf{x}\biggr)\\[-2pt]
&&\quad=P\biggl(-U\biggl(\frac1{T_1},\dots,\frac1{T_m}\biggr)\le
m\mathbf{x}\biggr)\\[-2pt]
&&\quad=\int_{S_m}P\biggl(-U\biggl(\frac1{t_1},\dots,\frac1{t_m}\biggr)\le
m\mathbf{x}\mid\bfT=\mathbf{t}\biggr)(P*\bfT)(\mathrm{d}\mathbf{t})\\[-2pt]
&&\quad=\int_{S_m}P\biggl(-U\biggl(\frac1{t_1},\dots,\frac1{t_m}\biggr)\le
m\mathbf{x}\biggr)\tilde\mu(\mathrm{d}\mathbf{t})\\[-2pt]
&&\quad=\frac1m \int_{S_m}P\biggl(-U\biggl(\frac1{t_1},\dots,\frac1{t_m}\biggr)\le
m\mathbf{x}\biggr)\mu(\mathrm{d}\mathbf{t})\\[-2pt]
&&\quad=\frac1m \int_{S_m}P\Bigl(U\ge m\max_{i\le m}(-x_it_i)\Bigr) \mu(\mathrm{d}\mathbf{t})\\[-2pt]
&&\quad=\frac1m \int_{S_m}P\biggl(U\ge-m\biggl(\sum_{j\le m}x_j\biggr)
\max_{i\le m}(\tilde x_it_i)\biggr) \mu(\mathrm{d}\mathbf{t})\\[-2pt]
&&\quad=\frac1m \int_{S_m} 1+ m\biggl(\sum_{j\le m}x_j\biggr) \max_{i\le
m}(\tilde x_it_i)\mu(\mathrm{d}\mathbf{t})\\[-2pt]
&&\quad=1+\biggl(\sum_{j\le m}x_j\biggr) \int_{S_m} \max_{i\le m}(\tilde
x_it_i)\mu(\mathrm{d}\mathbf{t}).
\end{eqnarray*}
This implies part (i) of the proposition.

On the other hand, we have for $\mathbf{x}\le\bfzero$ and large $s>0$
\begin{eqnarray*}
&&P\biggl(-U\biggl(\frac1{Z_1},\dots,\frac1{Z_m}\biggr)\le\frac1s\mathbf{x}\biggr)^s\\[-2pt]
&&\quad=\biggl(\int_{[\bfzero,\bfc]}P\biggl(U\ge\frac1s\max_{i\le
m}(-x_iz_i)\biggr)(P*\bfZ)(\mathrm{d}\mathbf{z})\biggr)^s\\[-2pt]
&&\quad=\biggl(1-\frac1s \int_{[\bfzero,\bfc]} \max_{i\le
m}(-x_iz_i)(P*\bfZ)(\mathrm{d}\mathbf{z})\biggr)^s\\[-2pt]
&&\quad\mathop{\to}\limits_{s\to\infty}\exp\biggl(-\int_{[\bfzero,\bfc]} \max_{i\le
m}(-x_iz_i)(P*\bfZ)(\mathrm{d}\mathbf{z})\biggr)\\[-2pt]
&&\quad=:G(\mathbf{x})
\end{eqnarray*}
with $\bfc=(c_1,\dots,c_m)$.\vadjust{\goodbreak}

Lemma 7.2.1 in Reiss \cite{re89} now implies that $G$ is a distribution
function that is obviously max stable: $G^s(s^{-1}\mathbf{x})=G(\mathbf{x})$,
$s>0$; that is, $G$ is a multivariate EVD and has negative standard
exponential margins $G_i(x)=\exp(xE(Z_i))=\exp(x)$, $x\le0$. As a
consequence, $1+\log(G(\mathbf{x}))$ is a GP function with
\begin{eqnarray*}
1+\log(G(\mathbf{x}))&=&1- \int_{[\bfzero,\bfc]} \max_{i\le
m}(-x_iz_i)(P*\bfZ)(\mathrm{d}\mathbf{z})\\
&=&P\biggl(-U\biggl(\frac1{Z_1},\dots,\frac1{Z_m}\biggr)\le\mathbf{x}\biggr)
\end{eqnarray*}
for $\mathbf{x}_0\le\mathbf{x}\le\bfzero$ and some $\mathbf{x}_0<\bfzero$.
\end{pf*}

Let, for instance, $C$ be an arbitrary $m$-dimensional \textit{copula};
that is, $C$ is the distribution function of a random vector $\bfS$
with uniform margins $P(S_i\le s)=s$, $s\in(0,1)$, $i\le m$, (Nelsen
\cite{nelsen06}). Then $\bfZ:=2\bfS$ is a proper choice in part (ii) of
Proposition~\ref{propgenerationofGPD}. Proposition~\ref{propgenerationofGPD},
therefore, maps the set of copulas in a
natural way to the set of multivariate GPDs, thus opening a wide range
of possible scenarios.

According to Theorem \ref{theocharacterizationofdomain}, we call a
copula $C_W$ a \textit{GPD copula on $[0,1]^m$} or simply a~\textit{GPD
copula} if there exists $\mathbf{y}_0<\mathbf{1}$ such that
\[
C_W(\mathbf{y}) = W(\mathbf{y}-\mathbf{1}),\qquad  \mathbf{y}_0\le\mathbf{y}\le\mathbf{1},
\]
where $W$ is a GPD with standard uniform margins in a left neighborhood
of zero.

For mathematical convenience we temporarily shift a copula to the
interval $[-1,0]^m$ by shifting each univariate margin by $-1$. Thus we
obtain a distribution function $\widetilde{C}_W$ from a GPD copula
$C_W$, whose marginal distribution functions are the uniform
distribution on $[-1,0]$, and $\widetilde{C}_W$
coincides close to zero with a GPD $W$ as in equation
\eqref{eqndehaanresnickpickandsrepresentationofgpd}; that is, there
exists $\mathbf{x}_0<\bfzero$ such that
\begin{eqnarray*}
\widetilde{C}_W(\mathbf{x})&=&W(x_1,\dots,x_m)\\
&=&1+\biggl(\sum_{j\le m}x_j\biggr)\int_{S_m} \max_{i\le m}\biggl(t_i
\frac{x_i}{\sum_{j\le m}x_j}\biggr)\mu(\mathrm{d}\mathbf{t}),\qquad \mathbf{x}\in[\mathbf
{x}_0,\bfzero].
\end{eqnarray*}
Because $\widetilde{C}_W$ inherits its properties from the original GPD
copula $C_W$, we call $\widetilde{C}_W$ a~\textit{GPD copula on
$[-1,0]^m$}.

For later purposes we remark that a random vector $\bfV\in[-1,0]^m$
following a GPD copula on $[-1,0]^m$ can easily be generated as
follows, using Proposition \ref{propgenerationofGPD}. Let $U$ be
uniformly distributed on $(0,1)$ and independent of the vector
$\bfS=(S_1,\dots,S_m)$, which follows an arbitrary copula on $[0,1]^m$.
Then we have for $i\le m$
\begin{eqnarray*}
P\biggl(-U\frac1{2S_i}\le x\biggr)&=&
\cases{
1+x,&\quad if
$-\dfrac12 \le x\le0$,\cr
\dfrac1{4\abs{x}},&\quad if $x<-\dfrac12$,
}\\
&=&\!:H(x),\qquad  x\le0,
\end{eqnarray*}
and, consequently,
\[
\bfV:=\biggl(H\biggl(-\frac{U}{2S_1}\biggr)-1,\dots,
H\biggl(-\frac{U}{2S_m}\biggr)-1\biggr)=(V_1,\dots,V_m)
\]
with
%
\begin{equation} \label{eqngenerationofGPD-copula}
V_i=
\cases{ -\dfrac{U}{2S_i},&\quad if $U\le S_i$,\cr
\dfrac{S_i}{2U}-1,&\quad if $U> S_i$,
}
\end{equation}
follows by Proposition \ref{propgenerationofGPD} a GPD copula on
$[-1,0]^m$.

\section{Multivariate piecing together}\label{secmultivariatepiecingtogether}
The multivariate PT approach consists of two steps. In a first step,
the upper tail of a given $m$-dimensional copula $C$ is cut off and
substituted by the upper tail of multivariate GPD copula in a
continuous manner. The result is again a copula, that is, an
$m$-dimensional distribution with uniform margins. The other step
consists of the transformation of each margin of this copula by a given
univariate distribution function $F_i^*$, $1\le i\le m$. This provides,
altogether, a multivariate distribution function with prescribed
margins $F_i^*$ whose copula coincides in its central part with $C$ and
in its upper tail with a GPD copula.

We start with fitting a GPD copula to the upper tail of a given copula
$C$ on $[-1,0]^m$. Recall that for mathematical convenience we shift
any copula $\tilde C(\bfu)$, $\bfu\in[0,1]^m$, to a~copula on
$[-1,0]^m$ by setting $C(\mathbf{v})=\tilde C(\bfone+\mathbf{v})$,
$\mathbf{v}\in[-1,0]^m$.

Let $\bfV=(V_1,\dots,V_m)$ follow a GPD copula on $[-1,0]^m$; that is,
$P(V_i\le x)=1+x$, $-1\le x\le0$, is for each $i\le m$ the uniform
distribution on $[-1,0]$, and there exists $\mathbf{x}_0=
(x_0^{(1)},\dots,x_0^{(m)} )<\bfzero$ such that for each
$\mathbf{x}=(x_1,\dots,x_m)\in[\mathbf{x}_0,\bfzero]$
\[
P(\bfV\le\mathbf{x})=1+\biggl(\sum_{i\le m}x_i\biggr) D\biggl(\frac{x_1}{\sum_{i\le
m}x_i},\dots,\frac{x_{m-1}}{\sum_{i\le m}x_i}\biggr),
\]
where $D$ is a Pickands dependence function.

Let $\bfY=(Y_1,\dots,Y_m)$ follow an arbitrary copula $C$ on $[-1,0]^m$
and suppose that $\bfY$ is independent of $\bfV$. Choose a threshold
$\mathbf{y}=(y_1,\dots,y_m)\in[-1,0]^m$ and put
%
\begin{equation}\label{eqndefinitionofpiecingtogether}
Q_i:=Y_i1_{(Y_i\le y_i)}-y_iV_i1_{(Y_i>y_i)},\qquad  i\le m.
\end{equation}
The random vector $\bfQ$ then follows a GPD copula on $[-1,0]^m$, which
coincides with $C$ on $\times_{i\le m}[-1,y_i]$. This is the content of
the main result of this section.

\begin{prop}
Suppose that $P(\bfY>\mathbf{y})>0$. Each $Q_i$ defined in
\eqref{eqndefinitionofpiecingtogether} follows the uniform distribution
on $[-1,0]$. The random vector $\bfQ=(Q_1,\dots,Q_m)$ follows a GPD
copula on $[-1,0]^m$, which coincides with $C$ on $\times_{i\le
m}[-1,y_i]$; that is,
\[
P(\bfQ\le\mathbf{x})=C(\mathbf{x}),\qquad \mathbf{x}\le\mathbf{y}.
\]
We have, moreover, with $x_i\in [\max (y_i,x_0^{(i)} ),0 ]$, $i\le m$,
for any non-empty subset $K$ of $\set{1,\dots,m}$
\[
P(Q_i\ge x_i,i\in K)=P(V_i\ge b_{i,K}x_i, i\in K),
\]
where
\[
b_{i,K}:=\frac{P(Y_j>y_j,j\in K)}{-y_i}= \frac{P(Y_j>y_j,j\in K
)}{P(Y_i>y_i)}\in(0,1],\qquad  i\in K.
\]
\end{prop}

\begin{pf}
First we show that each $Q_i$ follows the uniform distribution on
$[-1,0]$. We have for $-1\le x\le y_i$
\begin{eqnarray*}
P(Q_i\le x)&=&P(Q_i\le x,Y_i\le y_i)+P(Q_i\le x,Y_i>y_i)\\[-1pt]
&=&P(Y_i\le x)\\[-1pt]
&=&1+x,
\end{eqnarray*}
whereas for $y_i<x\le0$ we obtain
\begin{eqnarray*}
P(Q_i\le x)&=&P(Y_i\le y_i)+P(-y_iV_i\le x)P(Y_i>y_i)\\[-1pt]
&=&1+y_i+P\biggl(V_i\le-\frac{x}{y_i}\biggr)(-y_i)\\[-1pt]
&=&1+y_i+\biggl(1-\frac{x}{y_i}\biggr)(-y_i)\\[-1pt]
&=&1+x.
\end{eqnarray*}
The random vector $\bfQ$, thus, follows a copula on $[-1,0]^m$. We
have, further, for $\mathbf{x}\le\mathbf{y}$
\begin{eqnarray*}
P(\bfQ\le\mathbf{x})&=&P(\bfQ\le\mathbf{x},\bfY\le
\mathbf{y})+P(\bfQ\le\mathbf{x},\bfY\not\le\mathbf{y})\\[-1pt]
&=&P(\bfY\le\mathbf{x})\\[-1pt]
&=&C(\mathbf{x}).
\end{eqnarray*}
By Proposition 2.1 in Falk and Michel \cite{fami08} we have with
$x_i\in[\max(y_i,\omega_i),0]$, $i\le m$, $t\in[0,1]$ and an arbitrary
subset $K\subset\set{1,\dots,m}$
\begin{eqnarray*}
P(Q_j>tx_j,j\in K)&=&P(Q_j>tx_j,Y_j>y_j,j\in K)\\[-1pt]
&=&P(-y_jV_j>tx_j,j\in K)P(Y_j>y_j,j\in K)\\[-1pt]
&=&tP(-y_jV_j>x_j,j\in K)P(Y_j>y_j,j\in K)\\[-1pt]
&=&tP(Q_j>x_j,j\in K),
\end{eqnarray*}
which, again by Proposition 2.1 in Falk and Michel \cite{fami08}, implies that
$\bfQ$ follows a GPD.\vadjust{\goodbreak}

We have, moreover, with $x_i\in [\max (y_i,x_0^{(i)} ),0 ]$, $i\le m$,
\begin{eqnarray*}
&&P(Q_i\ge x_i,i\in K)\\
&&\quad=P(Q_i\ge x_i,Y_i>y_i, i\in K)+P(Q_i\ge x_i,i\in K,Y_j\le
y_j\mbox{ for some }j\in K)\\
&&\quad=P(-y_iV_i\ge x_i,i\in K)P(Y_i>y_i,i\in K)\\
&&\quad=P\biggl(V_i\ge-\frac{x_i}{y_i},i\in K\biggr)P(Y_i>y_i,i\in K)\\
&&\quad=P(V_i\ge b_{i,K}x_i,i\in K).
\end{eqnarray*}
\upqed\end{pf}

The above approach provides an easy way to generate a random vector
$\bfX\in\R^m$ with prescribed margins $F_i^*$, $i\le m$, such that
$\bfX$ has a given copula in the central part of the data, whereas in
the upper tail it has a GPD copula.

Take $\bfQ=(Q_1,\dots,Q_m)$ as in
\eqref{eqndefinitionofpiecingtogether} and put $\widetilde
\bfQ:=(Q_1+1,\dots,Q_m+1)$. Then each component $\widetilde Q_i$ of
$\widetilde\bfQ$ is uniformly distributed on $(0,1)$ and thus
%
\begin{equation}\label{eqnGPD-rvfollowingaprescribeddfbelowthreshold}
\bfX:=(X_1,\dots,X_m):=(F_1^{*-1}(\widetilde Q_1),\dots,
F_m^{*-1}(\widetilde Q_m))
\end{equation}
has the desired properties.

Combining the univariate \textit{and} the multivariate PT approach now
consists in choosing a threshold $u(i)\in\R$ for each dimension $i\le
m$ and a univariate distribution function $F_i$ together with an
arbitrary univariate GPD $W_{\gamma_i,\mu_i\sigma_i}$, and putting for
$i\le m$
%
\begin{equation}\label{eqndefinitionofmultivariatepiecingtogether}
F_i^*(x):=
\cases{ F_i(x),&\quad if $x\le u(i)$\cr
\bigl(1-F_i(u(i))\bigr)W_{\gamma_i,\mu_i\sigma_i}(x)+F_i(u(i)),&\quad if $x>u(i)$}.
\end{equation}
This is typically done in a way such that $F_i^*$ is a continuous
function.

\section{An application to operational loss data}\label{secoperationalrisk}
In this section we apply our
multivariate PT approach to operational loss data. For an excellent
introduction to operational risk and insurance analytics, see Chapter
10 of McNeil \textit{et al.} \cite{mcfrem05} and the literature cited there. In the
sequel we give a brief summary.

According to the New Basel Capital Accord (Basel II), banks are
required to determine the regulatory capital charge for operational
risk, defined as the risk of losses resulting from inadequate or failed
internal processes, people and systems or from external events. The
Basel Committee on Banking Supervision encourages the use and further
development of advanced modelling techniques to quantify operational
risk. The most risk-sensitive methodology is the \textit{loss
distribution approach} using bank internal data to estimate probability
distribution functions for each business line/event type category. To
provide a greater consistency of loss data collection within and
between banks, operational losses are classified in eight business
lines and seven event types. To calculate the capital charge for each
business line/event type combination, a risk measure such as value at
risk to the 99.9\% confidence level over a one-year holding period is
chosen. A~conservative way to assess a bank's total capital requirement
is to sum up the capital charges across business line/event type
classes assuming perfect dependence and disregarding diversification
effects in operational risk. Actually, the fact that all severe losses
occur in the same year is rather dubious. Therefore, the dependence
structure among losses of different business line/event type categories
needs to be modelled explicitly. For simplicity, we consider in what
follows only business lines and not event types.

The frequency of a loss event for business line $i$ over a one-year
time horizon will be denoted by $N(i)$. The random loss associated with
the $k$th loss event for business line $i$ will be denoted by
$\zeta_k(i)$.

The random loss $L(i)$ over one year for business line $i$ is,
therefore, modelled as
\[
L(i)=\sum_{k=1}^{N(i)}\zeta_k(i),
\]
where $\zeta_1(i), \zeta_2(i),\dots$ are assumed to be i.i.d. with
distribution function $F_i$ and they are independent of their total
number $N(i)$.

The goal is to model the total loss distribution for operational risk;
that is, the distribution of
\[
L:=\sum_{i=1}^m L(i),
\]
or parameters of it such as the value at risk $\operatorname{VAR}(\alpha)$ at
the probability level $\alpha$ satisfying
\[
P\bigl(L\ge\operatorname{VAR}(\alpha)\bigr)=1-\alpha
\]
or the expected shortfall at the probability level $\alpha$
\[
\operatorname{ES}(\alpha):=E\bigl(L\mid L\ge\operatorname{VAR}(\alpha)\bigr).
\]

In order to assess the total capital charge, the traditional models for
measuring operational risk determine $\operatorname{VAR}(\alpha)$ and
$\operatorname{ES}(\alpha)$ for each of the $m$ business lines separately and
then simply sum up the corresponding capital charges.

In contrast, Di Clemente and Romano \cite{diclmenenteromano} suggest modelling the
dependence structure among $L(1),\dots,L(m)$ by a copula function,
precisely, by the copula corresponding to the $m$-dimensional
$t$-distribution with $\nu$ degrees of freedom. For a closer look at
the issue of modelling the dependence among components of a random
vector of financial risk factors using the concept of a copula, see
Chapter 5 of McNeil \textit{et al.} \cite{mcfrem05}.

In our application we analyse operational losses of the external
database \textit{SAS OpRisk Global Data}, which contains worldwide
information on publicly reported operational losses over US \$$100\,000$.
Since we do not know the probability of losses lying under US
\$$100\,000$, we neglect this cut-off limit in modelling the severity and
frequency of the data. We concentrate on two business lines of the
financial sector, \textit{Commercial banking} and \textit{Retail
banking}.

First we follow the copula extreme value theory approach for modelling
operational loss data as outlined in Di Clemente and Romano \cite{diclmenenteromano}, but
we add the multivariate PT approach developed in Section
\ref{secmultivariatepiecingtogether}.

\subsection{Estimation}\label{secestimation}
Before applying the multivariate PT approach, we explore the
characteristics of the empirical distributions of the two business
lines and estimate the model's parameters. Thereby the severity
distribution of the random variable $\zeta_k(i)$ and the frequency
distribution of the random variable $N(i)$, $i = 1,2$, are treated
separately. To obtain the distribution function of the total loss
$L(i)$ of the business line~$i$ over a one-year time horizon
we accomplish a Monte Carlo
simulation combining the severity distribution with the frequency
distribution.

First we analyse the empirical distributions of the loss severity. The
measures skewness and kurtosis indicate that the empirical
distributions of the two business lines are skewed to the right and
very heavy tailed.

In the next step, parametric distributions (i.e., Weibull, gamma and
lognormal distribution) are fitted to the data. The parameters are
estimated by the maximum likelihood method.
%
%
%
With the help of graphical analysis (QQ plots, theoretical versus
empirical distribution function plots) and goodness-of-fit tests
(Anderson--Darling test, Cramer--von Mises (CvM) test), we conclude,
that none of the selected distributions provides a good fit to the
complete data sets. (For a detailed presentation and discussion of
goodness-of-fit techniques, see D'Agostino and Stephens \cite{daste86}.)
However, the lognormal distribution fits the body of the data very
well, while it underestimates the severity of the data in the right
tail.\looseness=-1

Therefore, to fit the tail data accurately, the univariate POT method
in the model of the severity distribution of the losses $\zeta_k(i)$ is
applied: The existence of a threshold~$u(i)$ for each business line $i$
is assumed such that $\zeta_k(i)$ follows a lognormal distribution
function below $u(i)$, whereas above $u(i)$ it follows a univariate
GPD, that is,
\begin{eqnarray} \label{eqnsimulationpiecingtogether}
&&P\bigl(\zeta_k(i)\le x\bigr)\nonumber\\ [-8pt]\\ [-8pt]
&&\quad =
\cases{
F_i(x),&\quad $x\le u(i)$,\cr
F_i(u(i)) +  \bigl(1-F_i(u(i)) \bigr) \operatorname{GPD}_{\beta(i),\xi(i)}\bigl(x-u(i)\bigr),&\quad $x\ge
u(i)$,}\nonumber
\end{eqnarray}
where the GPD is given by
\[
\operatorname{GPD}_{\beta(i),\xi(i)}(z):=
1-\biggl(1+\xi(i)\frac{z}{\beta(i)}\biggr)^{-1/\xi(i)},\qquad  z\ge0,
\]
with shape and scale parameter $\xi(i)>0$, $\beta(i)>0$. Furthermore,
$F_i$ is defined as $F_i(x) := \Phi ( (\log(x)-\mu(i) )/\sigma(i) )$,
where $\Phi$ is the standard normal distribution function and
$\mu(i)\in\R$, $\sigma(i)>0$ are location and scale parameters of the
lognormal distribution.

In this case we obtain for $x\ge u(i)$
\[
P\bigl(\zeta_k(i)>x\bigr)=P\bigl(\zeta_k(i)>u(i)\bigr)\biggl(1+\xi(i)\frac{x-u(i)}{\beta
(i)}\biggr)^{-1/\xi(i)}.\vadjust{\goodbreak}
\]

The threshold $u(i)$ is chosen with the help of mean excess plots. The
shape and scale parameters of the GPD are estimated by the maximum
likelihood method. For a discussion of the parameter estimation of a
GPD and optimal choice of the threshold, see Section 6 of Embrechts \textit{et al.} \cite{embrechts08}.

To determine the frequency distribution of the random variable $N(i)$,
the Poisson and negative binomial distribution are fitted to the total
number of losses per year. The parameters of these distributions are
estimated by the method of moments. Since the negative binomial
distribution has two parameters, $\alpha$ and $r$, it is more flexible
and often provides a better fit to operational loss data than the
Poisson distribution; see Cruz~\cite{cruz02}, page 89. With the help of the
$\chi^2$ goodness-of-fit test, this expectation is confirmed.
Therefore, the random variable $N(i)$ is modelled by the negative
binomial distribution, whose probability mass function is expressed as
\[
P\bigl(N(i)=n\bigr)= {\alpha(i)+n-1 \choose n} \biggl(\frac{1}{1+r(i)}\biggr)^{\alpha(i)}
\biggl(\frac{r(i)}{1+r(i)}\biggr)^n,\qquad   n\in\N_{0},
\]
with $\alpha(i)>0$, $r(i)>0$. The resulting estimates of the model's
parameters are given in Table \ref{tablemodelpara}.

%
\begin{table}
\tablewidth=350pt
\caption{Estimated model parameters} \label{tablemodelpara}
\begin{tabular*}{350pt}{@{\extracolsep{\fill}}llllllll@{}}
\hline
&$\alpha(i)$ & $r(i)$ & $\mu(i)$ & $\sigma(i)$ & $u(i)$ & $\beta(i)$ & $\xi (i)$\\
\hline
Commercial banking & 0.74 & \phantom{0}46.10 & 2.19 & 2.23 & 918.02 & 609.84 &
0.82\\
Retail banking & 0.39 & 162.04 & 0.88 & 2.06 & \phantom{0}69.18 & \phantom{0}99.75 & 1.02 \\
\hline
\end{tabular*}
\end{table}
%

\begin{table}
\tablewidth=230pt
\caption{Estimated correlation matrix for the $t$ copula}\label{tablecorrmatrix}
\begin{tabular*}{230pt}{@{\extracolsep{\fill}}lll@{}}
\hline
& Commercial & Retail \\
\hline
Commercial & 1 & 0.76 \\
Retail & 0.76 & 1 \\
\hline
\end{tabular*}
\end{table}

In the following, we model the dependence structure among $L(1)$ and
$L(2)$ by a copula function. The assumption of a normal copula has been
quite popular in finance for modelling the dependence between different
risks, but it puts less weight on observations that are large in each
component; see, for example, Rachev \textit{et al.} \cite{rachevsteinwei2009}. The $t$ copula is
more heavily tailed and, therefore, better suited for modelling
operational risk.

In our bivariate case the $t$ copula is fitted to the total loss data
over a one-year time horizon. The parameter of the correlation matrix
and the degrees of freedom $\nu$ are estimated by the maximum
likelihood method. For a discussion of the problem of fitting copulas
to data, see Section 5.5 of McNeil \textit{et al.} \cite{mcfrem05}.

Table \ref{tablecorrmatrix} contains the estimated correlation matrix for the
$t$ copula with $\nu=8.64$ estimated degrees of freedom.

To evaluate the goodness of fit of the $t$ copula, the CvM test is
applied. For recent reviews of copula goodness-of-fit testing, see Berg
\cite{berg09} and Genest \textit{et al.} \cite{gerebe09}. In Table \ref{tablecmtest} the CvM
test value and corresponding $p$-value are reported.

\begin{table}
\tablewidth=200pt
\caption{Goodness-of-fit test for the $t$ copula}\label{tablecmtest}
\begin{tabular*}{200pt}{@{\extracolsep{\fill}}ll@{}}
\hline
CvM statistic & $p$-value \\
\hline
0.02543851 & 0.521978\\
\hline
\end{tabular*}
\end{table}

Since the $p$-value is 0.521978, the null hypothesis that the
dependence structure of the data follows a $t$ copula is not rejected.


\subsection{Simulation}\label{secsimulation}
In the previous section we estimated the parameters that specify a
situation where only the univariate PT approach is applied. This model
is similar to the one described in Di Clemente and Romano \cite{diclmenenteromano}. Now
we show by simulations how popular risk measures such as value at risk
or expected shortfall can be influenced by the combination of
univariate and multivariate piecing together.

First we take the $t$ copula derived in Section \ref{secestimation}and
add the multivariate PT approach developed in Section
\ref{secmultivariatepiecingtogether}. We simulate $10^4$ independent
copies of $\widetilde{\bfY} = (\widetilde{Y}_1, \widetilde{Y}_2)$ and
$\bfV= (V_1, V_2),$ which follow the $t$ copula from above and a GPD
copula on $[-1,0]^2$, respectively. The realizations of $\bfY:=
\widetilde{\bfY} - \bfone$ and $\bfV$ are then combined with those
of a
random vector~$\bfQ$ according to definition
\eqref{eqndefinitionofpiecingtogether}. The distribution of $\bfQ$ is
then a GPD copula on $[-1,0]^2$ which coincides with the previously
mentioned $t$ copula -- shifted by minus one -- below some threshold
vector $\mathbf{y}= (y_1, y_2)$. The last step consists in shifting the
realizations of $\bfQ$ to the interval $[0,1]^2$ and transforming the
margins by $F_1^*, F_2^*$ according to equation
\eqref{eqnGPD-rvfollowingaprescribeddfbelowthreshold}. $(F_1^*,
F_2^*$ are derived from Monte Carlo simulations as described in Section
\ref{secestimation}.) Thus we obtain $10^4$ realizations of a random
vector $\bfX$ that follows a multivariate distribution function that
has the marginal distribution functions $F_1^*, F_2^*$ and the
associated copula is a GPD copula that coincides with the original $t$
copula in its central part. These realizations of $\bfX$ are then taken
to compute the empirical counterparts of the value at risk and the
expected shortfall.

Before we apply these steps, we remark that there are still two
remaining degrees of freedom in our model: the GPD copula that
underlies~$\bfV$ and the copula threshold vector~$\mathbf{y}$.
Our goal in this section is to give a first insight into the
consequences of~replacing the upper tail of a given copula with a GPD
copula. (Note that this procedure is justified by Theorem
\ref{theocharacterizationofdomain} and Corollary~\ref{coroequivalencesofdomain}.) For this purpose, we assume a simple
model:
\begin{longlist}
\item We define the GPD copula underlying $\bfV$ indirectly
by setting $\bfZ:= 2\bfS$ in
Proposition \ref{propgenerationofGPD}, where $\bfS$
follows a bivariate normal copula.
\item The copula threshold vector is obtained by $\mathbf{y}:=
 (F_1(u(1)), F_2(u(2)) )-\bfone$; that is, the
thresholds for the marginal distributions in the
univariate PT approach (see equation~\eqref{eqnsimulationpiecingtogether}) are
transformed and used for the multivariate PT
approach, too.
\end{longlist}

This way, we construct a parametric model that is only dependent on the
correlation matrix
\[
\bfSigma=
\pmatrix{
1 & \varrho\cr
\varrho& 1
},\qquad
 \varrho\in[0,1),
\]
of the random vector $\bfS$. It is typically very difficult,
particularly in higher dimensions, to find a good multivariate model
that describes both marginal behavior and dependence structure
effectively. The advantage of the preceding copula model is that it
depends on just one parameter $\varrho\in[0,1)$, allowing dependence
and independence of the margins in a~simple and continuous manner. We
refer again to McNeil \textit{et al.} \cite{mcfrem05}.

The simulations as described above were now done for various values of
$\varrho$. To attain more reliable estimates for the value at risk and
the expected shortfall, we simulated not only once but 50 times and
took the average. Additionally, these values were also computed for the
case in which only the univariate PT approach is applied and the
$t$ copula kept unchanged. 
This procedure allows us to identify the effect the multivariate PT
approach has on the mentioned risk measures.

In the following we state our results from the simulation series with
$\varrho=0.7,$ which models the case of relative high dependence but
not complete dependence between the business lines. Although a
graphical analysis of the used GPD copulas suggested that the degree of
dependence in the upper tail was increasing with $\varrho$ getting
larger, there was no observable trend in the estimates of the value at
risk and the
expected shortfall. 
Further research is necessary to derive criteria for the optimal choice
of $\varrho$ or, more generally, of the GPD copula underlying $\bfV$
and the copula threshold $\mathbf{y}$.

We now start with the presentation of the results. For simplicity, we
identify the business lines Commercial banking and Retail banking with
the cases $i=1$ and $i=2$, respectively. The random vector $\bfX$ from
above models the combined random losses for these two business lines
over one year, that is, $\bfX=  (L(1), L(2) )$. The value at risk and
the expected shortfall of $L(1)$, $L(2)$ and $L = L(1)+L(2)$ were
computed using their empirical counterparts
\[
\widehat{\operatorname{VAR}}(\alpha) = \widehat{F}^{-1}(\alpha),
\]
where $\widehat{F}$ is the empirical distribution function of $L(1)$,
$L(2)$ or $L$, respectively, and
\[
\widehat{\operatorname{ES}}(\alpha) = \frac{1}{n(1-\alpha)} \sum
_{i=1}^{n} l_i
1_{[\widehat{\operatorname{VAR}}(\alpha), \infty)}(l_i),
\]
where $l_i$ is the $i$th realization of $L(1)$, $L(2)$ or $L$,
respectively. By $1_B$ we denote the indicator function of a set $B$,
that is, $1_B(x)=1$ if $x\in B$ and $1_B(x)=0$ otherwise.

Table \ref{tableestimatedVARandESusingPTapproach} gives the means of
$50$ independent simulations resulting from the multivariate PT
approach, whereas Table \ref{tableestimatedVARandESusingunivPTapproach}
makes use of the univariate approach only.
%
\begin{table}
\caption{Estimated value at risk and expected shortfall,
GPD copula} \label{tableestimatedVARandESusingPTapproach}
\begin{tabular*}{\textwidth}{@{\extracolsep{\fill}}lllllllll@{}}
\hline
& \multicolumn{4}{c}{$\widehat{\operatorname{VAR}}(\alpha)$} &
\multicolumn{4}{c@{}}{$\widehat{\operatorname{ES}}(\alpha)$}\\ [-5pt]
& \multicolumn{4}{c}{\hrulefill} &
\multicolumn{4}{c@{}}{\hrulefill}\\
$\alpha$ & 95\% & 99\% & 99.5\% & 99.9\% & 95\% & 99\% & 99.5\% &
99.9\% \\
\hline
$L(1)$ & $13\,638$ & $32\,899$ & \phantom{0}$49\,650$ & $159\,442$ & \phantom{0}$46\,436$ & $154\,758$ &
$270\,036$ & $1038\,077$ \\
$L(2)$ & $12\,586$ & $45\,370$ & \phantom{0}$84\,386$ & $390\,127$ & \phantom{0}$93\,365$ & $381\,142$ &
$702\,350$ & $2880\,672$\\
$L$ & $26\,578$ & $75\,518$ & $127\,042$ & $533\,701$ & $135\,581$ & $512\,781$ & $930\,472$
& $3746\,889$ \\
\hline
\end{tabular*}
\end{table}
%
%
\begin{table}
\caption{Estimated value at risk and expected shortfall,
$t$ copula} \label{tableestimatedVARandESusingunivPTapproach}
\begin{tabular*}{\textwidth}{@{\extracolsep{\fill}}lllllllll@{}}
\hline
& \multicolumn{4}{c}{$\widehat{\operatorname{VAR}}(\alpha)$} &
\multicolumn{4}{c@{}}{$\widehat{\operatorname{ES}}(\alpha)$}\\ [-5pt]
& \multicolumn{4}{c}{\hrulefill} &
\multicolumn{4}{c@{}}{\hrulefill}\\
$\alpha$ & 95\% & 99\% & 99.5\% & 99.9\% & 95\% & 99\% & 99.5\% &
99.9\% \\
\hline
$L(1)$ & $13\,638$ & $32\,667$ & \phantom{0}$49\,196$ & $153\,322$ & \phantom{0}$37\,674$ & $111\,075$ &
$182\,956$ & \phantom{0$\,$}$608\,755$ \\ 
$L(2)$ & $12\,590$ & $45\,601$ & \phantom{0}$83\,414$ & $392\,673$ & \phantom{0}$71\,288$ & $270\,821$ &
$482\,058$ & $1774\,252$\\ 
$L$ & $25\,428$ & $75\,674$ & $131\,267$ & $533\,710$ & $105\,122$ & $366\,904$ & $637\,184$
& $2261\,144$ \\ 
\hline
\end{tabular*}
\end{table}
Since the marginal distributions are the same in both cases, namely
$F_1^*, F_2^*$, the value at risk estimates concerning
$L(1), L(2)$ nearly coincide across both tables. 
It is apparent that the respective values for $L$, the total loss, are
only slightly different.

On the other hand, the empirical expected shortfalls attain clearly
higher values if the upper tail of the $t$\vadjust{\goodbreak} copula is substituted by a
GPD. This behavior is independent of the~$\alpha$ level and holds for
the losses in the single business lines as well as for the total loss
in our simulation. This is remarkable since the corresponding estimates
for the value at risk are in both cases nearly the same, indicating
that there are some extreme high losses that occur very
rarely.\label{textextremelosses}

Clearly, with underlying estimated GPD shape parameters $\xi(1)=0.82$
and $\xi(2)=1.02$, the theoretical expected shortfall exists only for
$i=1$ but not for $i=2$. The significant increases for
$\widehat{\operatorname{ES}}(\alpha)$ in line $L(1)$ and $L(2)$ in Table
\ref{tableestimatedVARandESusingPTapproach} should, therefore, only be
due to the high volatility of the empirical expected shortfall, whereas
the significant increase in line $L$ should be caused by the
substituted GPD copula as well. This example of real operational loss
data might be considered as a warning, not to underestimate the effects
of rare events that nevertheless might occur simultaneously.

In cases of no existing theoretical expected shortfall, Moscadelli
\cite{mosca04} suggests the risk measure median shortfall (MS) that is defined
regardless of the values of the shape parameter $\xi$. If the
respective distribution function is continuous, the median shortfall
has the representation
\[
\operatorname{MS}(\alpha) = \operatorname{VAR} \biggl(\frac{1 + \alpha}{2} \biggr),
\]
see Biagini and Ulmer \cite{biagini09}, pages 749--750. In addition to our previous
results, Table \ref{tableestimatedMS} states the estimates for the
median shortfall, which we computed as $\widehat{\operatorname{MS}}(\alpha) =
\widehat{\operatorname{VAR}} ((1 + \alpha)/2 )$.
%
\begin{table}
\caption{Estimated median shortfall}\label{tableestimatedMS} 
\vspace*{-3pt}
\begin{tabular*}{\textwidth}{@{\extracolsep{\fill}}lllllllll@{}}
\hline
& \multicolumn{4}{c}{$\widehat{\operatorname{MS}}(\alpha)$, $t$
copula} & \multicolumn{4}{c@{}}{$\widehat{\operatorname{MS}}(\alpha)$, GPD
copula}\\ [-5pt]
& \multicolumn{4}{c}{\hrulefill} & \multicolumn{4}{c@{}}{\hrulefill}\\
$\alpha$ & 95\% & 99\% & 99.5\% & 99.9\% & 95\% & 99\% & 99.5\% &
99.9\% \\
\hline
$L(1)$ & $19\,829$ & \phantom{0}$49\,196$ & \phantom{0}$78\,678$ & $243\,938$ & $19\,829$ & \phantom{0}$49\,650$ & \phantom{0}$80\,292$
& \phantom{0}$313\,246$ \\
$L(2)$ & $21\,494$ & \phantom{0}$83\,414$ & $162\,438$ & $793\,252$ & $21\,600$ &
\phantom{0}$84\,386$ & $162\,866$ & \phantom{0}$782\,938$ \\
$L$ & $40\,340$ & $131\,267$ & $234\,910$ & $962\,458$ & $42\,463$ & $127\,042$ & $229\,260$
& $1085\,283$ \\
\hline
\end{tabular*}
\vspace*{-6pt}
\end{table}
%

Unlike the expected shortfall, the median shortfall as a robust measure
is not as heavily influenced by extreme values in the upper tail.
Therefore, the median shortfall estimates for $L(1), L(2)$ are closer
to each other than the expected shortfall estimates if we compare the
$t$ copula case with the GPD copula case.
%
%
However, the values for the median shortfall of $L$ at the 99.9\%\
level do clearly differ, indicating a small number of extreme high
losses in the upper tail -- in accordance with our considerations on
the expected shortfall, see page \pageref{textextremelosses}.

\begin{figure}

\includegraphics{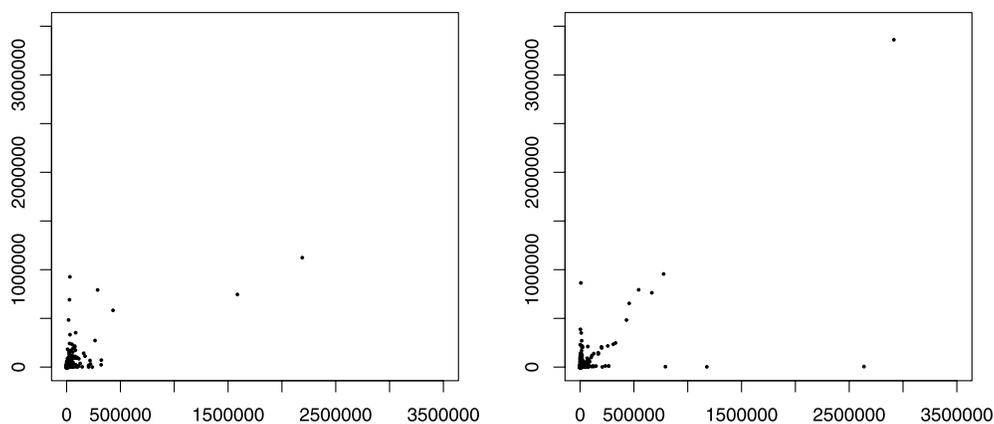}%
\vspace*{-3pt}
\caption{$10^4$ random deviates of $ (L(1), L(2) )$ based on the
original $t$ copula (left) and on the GPD copula (right).}
\label{figurecomparisondiClementeRomanoPT}
\end{figure}

To give visual insight into these results, Figure
\ref{figurecomparisondiClementeRomanoPT} compares the realizations of
$\bfX=  (L(1), L(2) )$ graphically.{}
%
The graphics were taken from a single simulation that contributed to
the results in Table \ref{tableestimatedVARandESusingPTapproach} and
Table \ref{tableestimatedVARandESusingunivPTapproach}. Although the $t$
copula itself is already heavily tailed, the substitution of its upper
part by a GPD puts even more weight on observations that are very high
in both components. (Recall Corollary
\ref{coroequivalencesofdomain}.) The latter type of modelling,
therefore, represents a higher risk of an extraordinarily high total
loss over a one-year time horizon. This can be seen on the point in the
upper right corner in the right scatterplot of Figure
\ref{figurecomparisondiClementeRomanoPT}, which represents a fictive
total loss of $6275\,000$, whereas the highest total loss in the
respective pure $t$ copula scenario is $3313\,000$.

\vspace*{-3pt} \section{Conclusions} \vspace*{-3pt}

In the present paper we extended the well known univariate PT approach
to higher dimensions. This was motivated by Theorem
\ref{theocharacterizationofdomain} and Corollary
\ref{coroequivalencesofdomain}, which show\vadjust{\goodbreak} that it is not sufficient to
apply the univariate approach to the marginal distributions of a random
vector if the upper tail of its distribution is to be modelled
adequately. It is, therefore, necessary to approximate the underlying
copula by a multivariate GPD.

The multivariate PT approach that was introduced in
\eqref{eqndefinitionofpiecingtogether} offers a wide range of scenarios
to be modelled because it depends basically only on some random vector
whose components need to be bounded and to have expectation one; see
Proposition \ref{propgenerationofGPD}. As a consequence we also
mentioned a natural way to map the set of copulas to the set of
multivariate GPDs that was useful for our simulation studies to obtain
a one-parametric model.

Because the values in a simulation are random, simulations occur that
produce no values that are high in both components. (This is depending
on the sample size, too.) Fixing this disadvantage is subject to
further research and could probably be achieved by making additional
restrictions on the GPD copula random vector in
\eqref{eqndefinitionofpiecingtogether}. Furthermore, goodness-of-fit
testing of the compound GPD copula is required.

Nevertheless, the multivariate PT approach is a powerful and suitable
tool to adequately model multivariate distributions in their upper
tails. This ensures that the probability of very rare events that occur
simultaneously and have a high effect if they occur is not
underestimated. The high empirical expected shortfalls in Table
\ref{tableestimatedVARandESusingPTapproach} might be considered as a
warning.

\section*{Acknowledgements}
The authors are grateful to the referees for their constructive
comments on an earlier version of this paper. The authors are also
indebted to the SAS Institute, Heidelberg, for providing access to the
database \textit{SAS OpRisk Global Data}.

The first author was supported by DFG Grant FA 262/4-1.

%

\printhistory


\begin{thebibliography}{35}

\bibitem{alnesch93}
\begin{barticle}[mr]
\bauthor{\bsnm{Alsina},~\bfnm{Claudi}\binits{C.}},
  \bauthor{\bsnm{Nelsen},~\bfnm{Roger~B.}\binits{R.B.}} \AND
  \bauthor{\bsnm{Schweizer},~\bfnm{Berthold}\binits{B.}}
(\byear{1993}).
\btitle{On the characterization of a~class of binary operations on distribution
  functions}.
\bjournal{Statist. Probab. Lett.}
\bvolume{17}
\bpages{85--89}.
\bid{doi={10.1016/0167-7152(93)90001-Y}, issn={0167-7152}, mr={1223530}}
\end{barticle}
\endbibitem

\bibitem{balkemahaan74}
\begin{barticle}[mr]
\bauthor{\bsnm{Balkema},~\bfnm{A.~A.}\binits{A.A.}} \AND
  \bauthor{\bparticle{de} \bsnm{Haan},~\bfnm{L.}\binits{L.}}
(\byear{1974}).
\btitle{Residual life time at great age}.
\bjournal{Ann. Probab.}
\bvolume{2}
\bpages{792--804}.
\bid{mr={0359049}}
\end{barticle}
\endbibitem

\bibitem{begosete04}
\begin{bbook}[mr]
\bauthor{\bsnm{Beirlant},~\bfnm{Jan}\binits{J.}},
  \bauthor{\bsnm{Goegebeur},~\bfnm{Yuri}\binits{Y.}},
  \bauthor{\bsnm{Teugels},~\bfnm{Jozef}\binits{J.}} \AND
  \bauthor{\bsnm{Segers},~\bfnm{Johan}\binits{J.}}
(\byear{2004}).
\btitle{Statistics of Extremes: Theory and Applications}.
\baddress{Chichester}: \bpublisher{Wiley}.
\bid{doi={10.1002/0470012382}, mr={2108013}}
\end{bbook}
\endbibitem


\bibitem{berg09}
\begin{bmisc}[author]
\bauthor{\bsnm{Berg},~\bfnm{D.}\binits{D.}}
(\byear{2009}).
\bhowpublished{Copula goodness-of-fit testing: An overview and
  power comparison. \textit{The European Journal of Finance}
  \textbf{15}
  675--701.}
\end{bmisc}
\endbibitem



\bibitem{biagini09}
\begin{bmisc}[author]
\bauthor{\bsnm{Biagini},~\bfnm{F.}\binits{F.}} \AND
\bauthor{\bsnm{Ulmer},~\bfnm{S.}\binits{S.}}
(\byear{2009}).
\bhowpublished{Asymptotics for operational
  risk quantified with expected shortfall. \textit{ASTIN Bulletin}
  \textbf{39} 735--752.}
\end{bmisc}
\endbibitem


\bibitem{buhazh06}
\begin{barticle}[mr]
\bauthor{\bsnm{Buishand},~\bfnm{T.~A.}\binits{T.A.}},
\bauthor{\bparticle{de} \bsnm{Haan},~\bfnm{L.}\binits{L.}} \AND
\bauthor{\bsnm{Zhou},~\bfnm{C.}\binits{C.}}
(\byear{2008}).
\btitle{On spatial extremes: With application to a rainfall problem}.
\bjournal{Ann. Appl. Statist.}
\bvolume{2}
\bpages{624--642}.
\bid{doi={10.1214/08-AOAS159}, issn={1932-6157}, mr={2524349}}
\end{barticle}
\endbibitem



\bibitem{cruz02}
\begin{bmisc}[author]
\bauthor{\bsnm{Cruz},~\bfnm{M.~G.}\binits{M.G.}}
(\byear{2002}).
\bhowpublished{\textit{Modeling, Measuring and Hedging
  Operational Risk}. Chichester: Wiley.}
\end{bmisc}
\endbibitem


\bibitem{daste86}
\begin{bbook}[mr]
\bauthor{\bsnm{D'Agostino},~\bfnm{L.}\binits{L.}} \AND
\bauthor{\bsnm{Stephens},~\bfnm{M.~A.}\binits{M.A.}}
(\byear{1986}).
\btitle{Goodness-of-fit Techniques}.
\bseries{Statistics: Textbooks and Monographs}
\bvolume{68}.
\baddress{New York}: \bpublisher{Marcel Dekker Inc.}
\bid{mr={0874534}}
\end{bbook}\vadjust{\goodbreak}
\endbibitem



\bibitem{haro98}
\begin{barticle}[mr]
\bauthor{\bparticle{de} \bsnm{Haan},~\bfnm{Laurens}\binits{L.}} \AND
  \bauthor{\bparticle{de} \bsnm{Ronde},~\bfnm{John}\binits{J.}}
(\byear{1998}).
\btitle{Sea and wind: Multivariate extremes at work}.
\bjournal{Extremes}
\bvolume{1}
\bpages{7--45}.
\bid{doi={10.1023/A:1009909800311}, issn={1386-1999}, mr={1652944}}
\end{barticle}
\endbibitem

\bibitem{hafe06}
\begin{bbook}[mr]
\bauthor{\bparticle{de} \bsnm{Haan},~\bfnm{Laurens}\binits{L.}} \AND
  \bauthor{\bsnm{Ferreira},~\bfnm{Ana}\binits{A.}}
(\byear{2006}).
\btitle{Extreme Value Theory: An Introduction}.
\baddress{New York}: \bpublisher{Springer}.
\bid{mr={2234156}}
\end{bbook}
\endbibitem

\bibitem{hasi99}
\begin{barticle}[mr]
\bauthor{\bparticle{de} \bsnm{Haan},~\bfnm{Laurens}\binits{L.}} \AND
  \bauthor{\bsnm{Sinha},~\bfnm{Ashoke~Kumar}\binits{A.K.}}
(\byear{1999}).
\btitle{Estimating the probability of a rare event}.
\bjournal{Ann. Statist.}
\bvolume{27}
\bpages{732--759}.
\bid{doi={10.1214/aos/1018031214}, issn={0090-5364}, mr={1714710}}
\end{barticle}
\endbibitem

\bibitem{deheuvels78}
\begin{bmisc}[author]
\bauthor{\bsnm{Deheuvels},~\bfnm{Paul}\binits{P.}}
(\byear{1978}).
\bhowpublished{Caract\'{e}risation compl\`{e}te des lois
  extr\^{e}me multivari\'{e}es et de la convergence des types extr\^{e}mes.
  \textit{Publ. Inst. Statist. Univ. Paris} \textbf{23} 1--36.}
\end{bmisc}
\endbibitem



\bibitem{deheuvels84}
\begin{bincollection}[mr]
\bauthor{\bsnm{Deheuvels},~\bfnm{Paul}\binits{P.}}
(\byear{1984}).
\btitle{Probabilistic aspects of multivariate extremes}.
In \bbooktitle{Statistical Extremes and Applications ({V}imeiro, 1983)}.
\bseries{NATO Adv. Sci. Inst. Ser. C Math. Phys. Sci.}
(\beditor{\bfnm{J. Tiago}\binits{J. Tiago} \bsnm{de Oliveira}}, ed.)
\bvolume{131}
\bpages{117--130}.
\baddress{Dordrecht}: \bpublisher{Reidel}.
\bid{mr={0784817}}
\end{bincollection}
\endbibitem

\bibitem{diclmenenteromano}
\begin{bmisc}[author]
\bauthor{\bparticle{Di} \bsnm{Clemente},~\bfnm{A.}\binits{A.}} \AND
  \bauthor{\bsnm{Romano},~\bfnm{C.}\binits{C.}}
  (\byear{2004}).
\bhowpublished{A copula-extreme value
  theory approach for modelling operational risk. In \textit{Operational Risk
  Modelling and Analysis, Theory and Practice} (M. Cruz, ed.) 189--208.
  London: RISK Books.}
\end{bmisc}
\endbibitem

\bibitem{embrechts08}
\begin{bbook}[mr]
\bauthor{\bsnm{Embrechts},~\bfnm{Paul}\binits{P.}},
  \bauthor{\bsnm{Kl{\"u}ppelberg},~\bfnm{Claudia}\binits{C.}} \AND
  \bauthor{\bsnm{Mikosch},~\bfnm{Thomas}\binits{T.}}
(\byear{1997}).
\btitle{Modelling Extremal Events for Insurance and Finance}.
\bseries{Applications of Mathematics (New York)}
\bvolume{33}.
\baddress{Berlin}: \bpublisher{Springer}.
\bid{mr={1458613}}
\bptok{imsref}%
\end{bbook}
\endbibitem


\bibitem{falk08}
\begin{barticle}[mr]
\bauthor{\bsnm{Falk},~\bfnm{Michael}\binits{M.}}
(\byear{2008}).
\btitle{It was 30 years ago today when {L}aurens de {H}aan went the
  multivariate way}.
\bjournal{Extremes}
\bvolume{11}
\bpages{55--80}.
\bid{doi={10.1007/s10687-007-0045-z}, issn={1386-1999}, mr={2420197}}
\end{barticle}
\endbibitem

\bibitem{fahure04}
\begin{bbook}[mr]
\bauthor{\bsnm{Falk},~\bfnm{Michael}\binits{M.}},
  \bauthor{\bsnm{H{\"u}sler},~\bfnm{J{\"u}rg}\binits{J.}} \AND
  \bauthor{\bsnm{Reiss},~\bfnm{Rolf-Dieter}\binits{R.D.}}
(\byear{2010}).
\btitle{Laws of Small Numbers: Extremes and Rare Events},
\bedition{3rd ed.}
\baddress{Basel}: \bpublisher{Birkh\"auser}.
\end{bbook}
\endbibitem

\bibitem{fami08}
\begin{barticle}[mr]
\bauthor{\bsnm{Falk},~\bfnm{M.}\binits{M.}} \AND
  \bauthor{\bsnm{Michel},~\bfnm{R.}\binits{R.}}
(\byear{2009}).
\btitle{Testing for a multivariate generalized {P}areto distribution}.
\bjournal{Extremes}
\bvolume{12}
\bpages{33--51}.
\bid{doi={10.1007/s10687-008-0067-1}, issn={1386-1999}, mr={2480722}}
\end{barticle}
\endbibitem

\bibitem{ga87}
\begin{bbook}[mr]
\bauthor{\bsnm{Galambos},~\bfnm{Janos}\binits{J.}}
(\byear{1987}).
\btitle{The Asymptotic Theory of Extreme Order Statistics},
\bedition{2nd ed.}
\baddress{Melbourne, FL}: \bpublisher{Robert E. Krieger Publishing Co. Inc.}
\bid{mr={0936631}}
\end{bbook}
\endbibitem

\bibitem{gemorose99}
\begin{barticle}[mr]
\bauthor{\bsnm{Genest},~\bfnm{C.}\binits{C.}},
  \bauthor{\bsnm{Quesada~Molina},~\bfnm{J.~J.}\binits{J.J.}},
  \bauthor{\bsnm{Rodr{\'{\i}}guez~Lallena},~\bfnm{J.~A.}\binits{J.A.}} \AND
  \bauthor{\bsnm{Sempi},~\bfnm{C.}\binits{C.}}
(\byear{1999}).
\btitle{A characterization of quasi-copulas}.
\bjournal{J. Multivariate Anal.}
\bvolume{69}
\bpages{193--205}.
\bid{doi={10.1006/jmva.1998.1809}, issn={0047-259X}, mr={1703371}}
\end{barticle}
\endbibitem

\bibitem{gerebe09}
\begin{barticle}[mr]
\bauthor{\bsnm{Genest},~\bfnm{Christian}\binits{C.}},
  \bauthor{\bsnm{R{\'e}millard},~\bfnm{Bruno}\binits{B.}} \AND
  \bauthor{\bsnm{Beaudoin},~\bfnm{David}\binits{D.}}
(\byear{2009}).
\btitle{Goodness-of-fit tests for copulas: A review and a power study}.
\bjournal{Insurance Math. Econom.}
\bvolume{44}
\bpages{199--213}.
\bid{doi={10.1016/j.insmatheco.2007.10.005}, issn={0167-6687}, mr={2517885}}
\end{barticle}
\endbibitem


\bibitem{huang92}
\begin{bmisc}[author]
\bauthor{\bsnm{Huang},~\bfnm{X.}\binits{X.}}
(\byear{1992}).
\bhowpublished{Statistics of bivariate extreme
  values. Ph.D. thesis, Tinbergen Institute Research Series.}
\end{bmisc}
\endbibitem



\bibitem{mcfrem05}
\begin{bbook}[mr]
\bauthor{\bsnm{McNeil},~\bfnm{Alexander~J.}\binits{A.J.}},
  \bauthor{\bsnm{Frey},~\bfnm{R{\"u}diger}\binits{R.}} \AND
  \bauthor{\bsnm{Embrechts},~\bfnm{Paul}\binits{P.}}
(\byear{2005}).
\btitle{Quantitative Risk Management: Concepts, Techniques and Tools}.
\baddress{Princeton, NJ}: \bpublisher{Princeton Univ. Press}.
\bid{mr={2175089}}
\end{bbook}
\endbibitem

\bibitem{mcneilneslehova09}
\begin{barticle}[mr]
\bauthor{\bsnm{McNeil},~\bfnm{Alexander~J.}\binits{A.J.}} \AND
  \bauthor{\bsnm{Ne{\v{s}}lehov{\'a}},~\bfnm{Johanna}\binits{J.}}
(\byear{2009}).
\btitle{Multivariate {A}rchimedean copulas, {$d$}-monotone functions and {$l\sb
  1$}-norm symmetric distributions}.
\bjournal{Ann. Statist.}
\bvolume{37}
\bpages{3059--3097}.
\bid{doi={10.1214/07-AOS556}, issn={0090-5364}, mr={2541455}}
\end{barticle}
\endbibitem

\bibitem{michel06}
\begin{bmisc}[author]
\bauthor{\bsnm{Michel},~\bfnm{Ren{\'e}}\binits{R.}}
(\byear{2006}).
\bhowpublished{Simulation and estimation in
  multivariate generalized Pareto models. Ph.D. thesis, Univ.
  W\"{u}rzburg.}
\end{bmisc}
\endbibitem


\bibitem{michel07}
\begin{barticle}[mr]
\bauthor{\bsnm{Michel},~\bfnm{Ren{\'e}}\binits{R.}}
(\byear{2007}).
\btitle{Simulation of certain multivariate generalized {P}areto distributions}.
\bjournal{Extremes}
\bvolume{10}
\bpages{83--107}.
\bid{doi={10.1007/s10687-007-0036-0}, issn={1386-1999}, mr={2415635}}
\end{barticle}
\endbibitem

\bibitem{michel08}
\begin{barticle}[mr]
\bauthor{\bsnm{Michel},~\bfnm{Ren{\'e}}\binits{R.}}
(\byear{2008}).
\btitle{Some notes on multivariate generalized {P}areto distributions}.
\bjournal{J. Multivariate Anal.}
\bvolume{99}
\bpages{1288--1301}.
\bid{doi={10.1016/j.jmva.2007.08.007}, issn={0047-259X}, mr={2432486}}
\end{barticle}
\endbibitem


\bibitem{mosca04}
\begin{bmisc}[author]
\bauthor{\bsnm{Moscadelli},~\bfnm{M.}\binits{M.}}
(\byear{2004}).
\bhowpublished{The modelling of operational risk:
  Experience with the analysis of the data collected by the Basel Committee.
  \textit{Banca D'Italia, Termini di discussione No.} \textbf{517}.}
\end{bmisc}
\endbibitem


\bibitem{nelsen06}
\begin{bbook}[mr]
\bauthor{\bsnm{Nelsen},~\bfnm{Roger~B.}\binits{R.B.}}
(\byear{2006}).
\btitle{An Introduction to Copulas},
\bedition{2nd ed.}
\baddress{New York}: \bpublisher{Springer}.
\bid{mr={2197664}}
\end{bbook}
\endbibitem

\bibitem{pickands75}
\begin{barticle}[mr]
\bauthor{\bsnm{Pickands},~\bfnm{James}\binits{J.} \bsuffix{III}}
(\byear{1975}).
\btitle{Statistical inference using extreme order statistics}.
\bjournal{Ann. Statist.}
\bvolume{3}
\bpages{119--131}.
\bid{issn={0090-5364}, mr={0423667}}
\end{barticle}
\endbibitem


\bibitem{rachevsteinwei2009}
\begin{bmisc}[author]
\bauthor{\bsnm{Rachev},~\bfnm{S.~T.}\binits{S.T.}},
\bauthor{\bsnm{Stein},~\bfnm{M.}\binits{M.}} \AND
\bauthor{\bsnm{Sun},~\bfnm{W.}\binits{W.}}
(\byear{2009}).
\bhowpublished{Copula concepts in
  financial markets. Technical report, Karslruhe Institute of Technology
  (KIT).}
\end{bmisc}
\endbibitem



\bibitem{re89}
\begin{bbook}[mr]
\bauthor{\bsnm{Reiss},~\bfnm{R.~D.}\binits{R.D.}}
(\byear{1989}).
\btitle{Approximate Distributions of Order Statistics:
With Applications to Nonparametric Statistics}.
\baddress{New York}: \bpublisher{Springer}.
\bid{mr={0988164}}
\end{bbook}
\endbibitem

\bibitem{reth}
\begin{bbook}[mr]
\bauthor{\bsnm{Reiss},~\bfnm{R.~D.}\binits{R.D.}} \AND
  \bauthor{\bsnm{Thomas},~\bfnm{M.}\binits{M.}}
(\byear{2007}).
\btitle{Statistical Analysis of Extreme Values},
\bedition{3rd ed.}
\baddress{Basel}: \bpublisher{Birkh\"auser}.
\bid{mr={1819648}}
\bptnote{check year}
\end{bbook}
\endbibitem

\bibitem{resnick06}
\begin{bbook}[mr]
\bauthor{\bsnm{Resnick},~\bfnm{Sidney~I.}\binits{S.I.}}
(\byear{2006}).
\btitle{Heavy-Tail Phenomena: Probabilistic and Statistical Modeling}.
\baddress{New York}: \bpublisher{Springer}.
\bid{mr={2271424}}
\bptnote{check year}
\end{bbook}
\endbibitem

\bibitem{rootzen05}
\begin{barticle}[mr]
\bauthor{\bsnm{Rootz{\'e}n},~\bfnm{Holger}\binits{H.}} \AND
  \bauthor{\bsnm{Tajvidi},~\bfnm{Nader}\binits{N.}}
(\byear{2006}).
\btitle{Multivariate generalized {P}areto distributions}.
\bjournal{Bernoulli}
\bvolume{12}
\bpages{917--930}.
\bid{doi={10.3150/bj/1161614952}, issn={1350-7265}, mr={2265668}}
\end{barticle}
\endbibitem

\end{thebibliography}
\end{document}